\documentclass[10pt]{amsart}

\newtheorem{theorem}{Theorem}
\newtheorem{lemma}{Lemma}

\newtheorem{cor}{Corollary}
\newtheorem{prop}{Proposition}

\usepackage{amsmath}
\usepackage[psamsfonts]{amssymb}
\usepackage[mathscr]{euscript}

\newcommand{\cal}{\mathcal}

\title{Pseudo-harmonic Hermitian structures on Weyl manifolds}
\author{Kamran Shakoor, Johann Davidov}

\address{Kamran Shakoor\\ Abdus Salam School of Mathematical Sciences, Government College University, Lahore, Pakistan.}
{\email{kamranshakoor@sms.edu.pk}

\address{Johann Davidov\\Institute of Mathematics and Informatics \\
Bulgarian Academy of Sciences\\ Acad. G. Bonchev St. Bl. 8\\ 1113 Sofia\\
Bulgaria }{\email{jtd@math.bas.bg}

\begin{document}

\begin{abstract}
We find geometric conditions on a Hermitian-Weyl manifold  under
which the  complex structure is a pseudo-harmonic map in the sense
of G. Kokarev \cite{K09} from the manifold into its twistor space.
This is done under the assumption that the dimension of the manifold
is four or the Hermitian-Weyl structure is locally conformally
K\"ahler.

\vspace{0,1cm} \noindent 2020 {\it Mathematics Subject Classification.} 53C43, 58E20.

\vspace{0,1cm} \noindent {\it Keywords: Almost-complex structures,
twistor spaces,  Weyl connections, pseudo-harmonic maps}

\end{abstract}

\thispagestyle{empty}

\maketitle

\section{Introduction}

Recall that an almost-complex structure on a Riemannian manifold is
called almost-Hermitian or compatible if it is an orthogonal
endomorphism of the tangent bundle of the manifold. It is well-known
that if a Riemannian manifold admits an almost-Hermitian structure,
it possesses many such structures, see, for example, \cite{D17}.
Thus, it is natural to look for \textquotedblleft reasonable"
criteria that distinguish some of these structures. One way to
obtain such criteria is to consider the almost-Hermitian structures
on a Riemannian manifold $(M,g)$ as sections of its twistor bundle
that parametrizes the almost-Hermitian structures on $(M,g)$. This
is the bundle $\pi:{\mathcal Z}\to M$ whose fibre at a point $p\in
M$ consists of all $g$-orthogonal complex structures $I_p:T_pM\to
T_pM$ ($I_p^2=-Id$) on the tangent space of $M$ at $p$. The fibre of
the bundle ${\mathcal Z}$ is the compact Hermitian symmetric space
$O(2m)/U(m)$, $2m=dim\,M$, and its standard metric
$G=-\frac{1}{2}Trace\,(I_1\circ I_2)$ is K\"ahler-Einstein. The
Levi-Civita connection of $(M,g)$ gives rise to a splitting
$T{\mathcal Z}={\mathcal H}\oplus {\mathcal V}$ of the tangent
bundle of ${\mathcal Z}$ into horizontal and vertical parts. This
decomposition allows one to define a $1$-parameter family of
Riemannian metrics $\widetilde g_t=\pi^{\ast}g+tG$, $t>0$, for which
the projection map $\pi:({\mathcal Z},\widetilde g_t)\to (M,g)$ is a
Riemannian submersion with totally geodesic fibres. In the
terminology of \cite[Definition 9.7]{Besse}, the family $\widetilde
g_t$ is the canonical variation of the metric $g$. Motivated by the
harmonic map theory, C. M. Wood \cite{W1,W2} has suggested to
consider as \textquotedblleft optimal" those almost-Hermitian
structures $J:(M,g)\to ({\mathcal Z },\widetilde g_1)$ that are
critical points of the energy functional under variations through
sections of ${\mathcal Z}$. In general, these critical points are
not harmonic maps, but, by analogy, they are referred to as
\textquotedblleft harmonic almost-complex structures" in \cite{W2};
they are also called \textquotedblleft harmonic sections" in
\cite{W1}, a term which seems more appropriate. Forgetting the
bundle structure of ${\mathcal Z}$, we can consider the
almost-Hermitian structures that are critical points of the energy
functional under variations through arbitrary maps $M\to{\mathcal
Z}$, not just sections. These structures are genuine harmonic maps
from $(M,g)$ into $({\mathcal Z}, \widetilde g_t)$, and we refer to
\cite{EL} for basic facts about harmonic maps. This point of view is
taken in \cite{DHM, Dav} where the problem of when an
almost-Hermitian structure on a Riemannian four-manifold is a
harmonic map from the manifold into its twistor space is discussed.
In \cite{DHM}, the metric $\widetilde g_t$ is defined via the
Levi-Civita connection; it is defined by means of a metric
connection with totally skew-symmetric torsion in \cite{Dav}. In
\cite{DM18} (cf. also \cite{D17}), the Riemannian $4$-manifolds
$(M,g)$ for which the Atiyah-Hitchin-Singer \cite{AHS} or the
Eells-Salamon   almost-complex structure \cite{ES} is a harmonic map
from the twistor space $({\mathcal Z},\widetilde g_t)$ of $(M,g)$
into the twistor space of $({\mathcal Z},\widetilde g_t)$ are
described.

If $g$ and $g^{1}=e^fg$ are conformal metrics, then clearly every
$g$-orthogonal endomorphism of the tangent bundle $TM$ is
$g^{1}$-orthogonal and vice versa. Hence the Riemannian manifolds
$(M,g)$ and $(M,g^{1})$ have the same twistor space, i.e. the
twistor space depends on the conformal class of $g$ rather than the
metric $g$ itself. Thus, it is natural to consider the twistor
spaces in the context of conformal geometry, see, for example,
\cite{AHS, Ga91}. The harmonic map techniques has a useful extension
in the conformal geometry introduced by G. Kokarev \cite{K09}.
Recall that a smooth map $\varphi:(M,g)\to (M',g')$ between
Riemannian manifolds is harmonic exactly when the trace of its
second fundamental form defined by means of the Levi-Civita
connections vanishes. If $M$ and $M'$ are endowed with torsion-free
connections $D$ and $D'$ respectively, a smooth map $\varphi:M\to
M'$ is called {\it pseudo-harmonic} in \cite{K09} if the trace of
its second fundamental form defined by means of $D$ and $D'$
vanishes. The pseudo-harmonic maps share some important properties
with the harmonic maps \cite{K09}, to mention here only the unique
continuation one.

Given a conformal  manifold $(M,\mathfrak{c})$, it is natural
to consider connections preserving the conformal class $\frak{c}$.
 Of special interest are those of them that have vanishing torsion.
 These are the Weyl connections studied in different aspects by many authors.
 Any Weyl connection $D$ gives rise to a splitting of the tangent bundle $T{\mathcal Z}$
 into horizontal and vertical parts. Then, taking a metric $g\in{\mathfrak c}$,
 one can define a $1$-parameter family $\widetilde g_t$ of Riemannian metrics on ${\mathcal Z}$ as above.
 We modify the Levi-Civita connection of $\widetilde g_t$ to
 a new torsion-free connection $D'$ on ${\mathcal Z}$, and find the conditions on $M$
 under which the map $J:M\to {\mathcal Z}$ defined by an
 integrable almost-Hermitian structure $J$ on $M$ is pseudo-harmonic
 with respect to $D$ and $D'$. This is done in the case when
 the dimension of $M$ is four, and  in higher dimension when the structure $(g,J)$ is locally conformally K\"ahler,
 i.e. when the Lee form $\theta$ of $(M,g,J)$ satisfies the identity $d\Omega=\theta\wedge\Omega$,
 where $\Omega(X,Y)=g(JX,Y)$ is the fundamental $2$-form of
 $(M,g,J)$. In dimension $4$, this identity is automatically
 satisfied. If  $dim\,M\geq 6$, it implies $d\theta=0$.  Note also
 that, in any dimension, a Hermitian structure is locally conformally
 K\"ahler if and only if $d\Omega=\theta\wedge\Omega$ and
 $d\theta=0$, see, for example, \cite{DO,V76}.
 If $D$ is the Weyl connection determined by $g$ and $\theta$, the map $J$ is always pseudo-harmonic.
 Two examples illustrating the obtained result are given,
one of them showing that there can be many Weyl connections for
which $J$ is a pseudo-harmonic map.

\section{Basics about twistor spaces}

In this section, we recall some basic facts about twistor spaces.
For more details, see, for example, \cite{D05, D17}.

\subsection{The manifold of compatible linear complex structures}\label{CLCS}

Let $V$ be a real vector space of even dimension $n=2m$ endowed with
an Euclidean metric $g$. Denote by $F(V)$ the set of all complex
structures on $V$ compatible with the metric $g$, i.e.
$g$-orthogonal. This set has the structure of an imbedded
submanifold of the vector space $\mathfrak{so}(V)$ of skew-symmetric
endomorphisms of $(V,g)$. The tangent space of $F(V)$ at a point $J$
consists of all endomorphisms $Q\in \mathfrak{so}(V)$ anti-commuting
with $J$, and we have the decomposition
$\mathfrak{so}(V)=T_JF(V)\oplus \{S\in\mathfrak{so}(V):~ SJ-JS=0\} $
which is orthogonal with respect to the metric
$G(S,T)=-\frac{1}{2}Trace\,(S\circ T)$ of $\mathfrak{so}(V)$.

The group $O(V)\cong O(2m)$ of orthogonal transformations of $V$ acts on $F(V)$
by conjugation and the isotropy subgroup at a fixed complex structure $J_0$ is isomorphic to the unitary group $U(m)$.
Therefore $F(V)$ can be identified with the homogeneous space $O(2m)/U(m)$. In particular, $dim\,F(V)=m^2-m$.
Note also that the manifold $F(V)$ has two connected components $F_{\pm}(V)$: if we fix an orientation on $V$,
these components consists of all complex structures on $V$ compatible with the metric $g$ and inducing $\pm$
the orientation of $V$; each of them has the homogeneous
representation $SO(2m)/U(m)$.

\smallskip

The metric $g$ of $V$ induces a metric on $\Lambda^2V$ given by
\begin{equation}\label{g on Lambda2}
g(v_1\wedge v_2,v_3\wedge
v_4)=\frac{1}{2}[g(v_1,v_3)g(v_2,v_4)-g(v_1,v_4)g(v_2,v_3)].
\end{equation}
Then we have an isomorphisms $\mathfrak{so}(V)\cong \Lambda^2V$, which sends $S\in \mathfrak{so}(V)$ to the $2$-vector $S^{\wedge}$ for which
$$
2g(S^{\wedge},u\wedge v)=g(S u,v),\quad u,v\in V.
$$
This isomorphism is an isometry with respect to the metric $\frac{1}{2}G$ on $\mathfrak{so}(V)$ and the metric $g$ on $\Lambda^2V$.

\subsection {The twistor space of a conformal manifold}

The twistor space of a conformal manifold $(M, {\mathfrak c})$ of even dimension $n=2m$
is the bundle $\pi:{\mathcal Z}\to M$ whose fibre at a point $p\in M$
consists of all complex structures of the tangent space $T_pM$ that are compatible
with the metrics of the conformal class ${\mathfrak c}$. If $M$ is oriented, we can consider the bundles over $M$
whose fibre at $p$ is the manifold of compatible complex structures yielding the positive and, respectively,
the negative orientation of $T_pM$. The latter bundles are disjoint open subsets of ${\mathcal Z}$ and are frequently called the positive
and the negative twistor spaces of $(M,{\mathfrak c})$, respectively. If $M$ is connected and oriented, they are the connected components of the manifold ${\mathcal Z}$.

The bundle ${\mathcal Z}$ can be considered as a subbundle of the
bundle  $A(TM)$ of the endomorphisms of $TM$ that are skew-symmetric
with respect to one Riemannian metric in ${\mathfrak c}$, hence with
respect to all such metrics. Let $D$ be a connection on
$(M,{\mathfrak c})$ preserving the conformal class in the sense that
for every Riemannian metric $g$ in the conformal class ${\mathfrak
c}$, there exists a $1$-form $\varphi_g$ such that
$Dg=\varphi_g\otimes g$; obviously, such a form is unique. The
connection on the vector bundle $Hom (TM,TM)$ induced by $D$ will
also be denoted by $D$. Endow the bundle $Hom(TM,TM)$ with the
metric $ G(a,b)=\frac{1}{2}Trace_g\{TM\ni X\to g(aX,bX)\}$, where
$g$ is a Riemannian metric in ${\mathfrak c}$. Clearly the metric
$G$ does not depend on the particular choice of $g$.

\begin{lemma}\label{D on Hom}
The induced connection $D$ preserves the metric $G$ and the bundle $A(TM)$.
\end{lemma}

\begin{proof}
Fix a point $p\in M$, and take a metric $g\in {\mathfrak c}$. Let
$e_1,...,e_n$ be a $g$-orthonormal basis of $T_pM$. We can extend
this basis to a frame of vector fields $E_1,...,E_{n}$ in a
(geodesically convex) neighbourhood of the point $p$ such that
$(E_i)_p=e_i$, $D E_i|_p=0$ and $g(E_i,E_j)=f\delta_{ij}$,
$i,j=1,...,n$,  where $f=f_{g}$ is a positive smooth function
depending on $g$. Let $a$, $b$ be sections of $Hom(TM,TM)$ and $X\in
T_pM$. Setting $\varphi=\varphi_{g}$, we have
$$
\begin{array}{l}
\displaystyle{(D_{X}G)(a,b)=\frac{1}{2}\sum\limits_{i=1}^n\big[X(\frac{1}{f}g(aE_i,bE_i))}
\displaystyle{\hfill-\frac{1}{f(p)}g(D_{X}(aE_i),bE_i)-\frac{1}{f(p)}g(aE_i,D_{X}(bE_i))\big]}\\[6pt]
\displaystyle{=X(\frac{1}{f})G(a,b)_p+\frac{1}{2}\frac{1}{f(p)}\sum\limits_{i=1}^n\big[X(g(aE_i,bE_i))}
\displaystyle{\hfill-g(D_{X}(aE_i),bE_i)-g(aE_i,D_{X}(bE_i))\big]},\\[6pt]
\displaystyle{=[X(\frac{1}{f})+\frac{1}{f(p)}\varphi(X)]G(a,b)_p}.
\end{array}
$$
Clearly, $f(p)=g(E_i,E_i)_p=g(e_i,e_i)=1$ and $X(f)=X(g(E_i,E_i))=\varphi(X)g(e_i,e_i)=\varphi(X)$. Therefore, $(D_{X}G)(a,b)=0$.

Let $S$ be a section of $A(TM)$ and $X,Y,Z$ vector fields on $M$. Then
$$
\begin{array}{c}
0=X(g(SY,Z)+g(Y,SZ))\\[6pt]
=\varphi(X)g(SY,Z)+g(D_{X}SY,Z)+g(SY,D_{X}Z)\\[6pt]
+\varphi(X)g(Y,SZ)+g(D_{X}Y,SZ)+g(Y,D_{X}SZ)\\[6pt]
=g((D_{X}S)Y,Z)+g(Y,(D_{X}S)Z).
\end{array}
$$
Thus, the endomorphism $D_{X}S$ of $TM$ is ${\mathfrak c}$-skew-symmetric.
\end{proof}

\section{Basics about Weyl connections}

Let $M$ be a manifold of dimension $n=2m$ endowed with a conformal
class ${\mathfrak c}$ of Riemannian metrics. Recall that  a Weyl
connection on $(M,{\mathfrak c})$  is a torsion-free connection $D$
on $M$ that preserves the conformal structure ${\mathfrak c}$. A
conformal manifold equipped with a Weyl connection is called a Weyl
manifold.

For $g\in {\mathfrak c}$, denote the Levi-Civita connection of $g$ by $\nabla^g$. Then
\begin{equation}\label{Weyl}
D_{X}Y=\nabla^g_{X}Y-\frac{1}{2}[\varphi_g(X)Y+\varphi_g(Y)X-g(X,Y)\varphi_g^{\sharp}],
\end{equation}
where the $1$-form $\varphi_g$ is determined by $Dg=\varphi_g\otimes
g$ and $\varphi_g^{\sharp}$ is the dual vector field of the form
$\varphi_g$ with respect to the metric $g$, i.e.
$g(\varphi_g^{\sharp},Z)=\varphi_g(Z)$. If $g^1=e^{f}g$ is another
Riemannian metric in the conformal class, $f$ being a smooth
function, then
\begin{equation}\label{change}
\begin{array}{c}
\nabla^{g^1}_{X}Y=\nabla^g_{X}Y+\displaystyle{\frac{1}{2}}[X(f)Y+Y(f)X-g(X,Y)\nabla^gf],\\[6pt]
\varphi_{g^1}=df+\varphi_{g},\quad\varphi_{g^1}^{\sharp}=e^{-f}(\varphi_{g}^{\sharp}+\nabla^g f),
\end{array}
\end{equation}
where $\nabla^g f$ is the gradient of $f$ with respect to $g$. In particular,
the condition $d\varphi_g=0$ does not depend on the choice of the metric $g$ in ${\mathfrak c}$.
In this case, we say that $D$ determines a closed Weyl structure. If $d\varphi_g=0$,
then locally $\varphi_g=d\psi$ for a smooth function $\psi$, so $D$ coincides locally with the Levi-Civita connection of the metric $e^{-\psi}g$.
The condition that the form $\varphi_g$ is exact also does not depend on the choice of $g$, and a Weyl structure with exact $\varphi_g$ is called exact.

If $g$ is a Riemannian metric and $\varphi$ is a $1$-form on a manifold, there is a unique Weyl connection for the conformal class ${\mathfrak c}$ of $g$
such that $\varphi_{g}=\varphi$. Indeed, for any conformal metric $g^1=e^{f}g$, set $\varphi_{g^1}=df+\varphi_{g}$.
Then, by (\ref{change}), the right-hand side of (\ref{Weyl}) does not depend on the choice of a metric in ${\mathfrak c}$,
so identity (\ref{Weyl}) defines a Weyl connection.

\smallskip

Now, let $D$ be a Weyl connection on $(M,{\mathfrak c})$. Let $R^D$
and $R^g$ be the curvature tensors of the connections $D$ and
$\nabla^g$, $g\in {\mathfrak c}$, respectively.

\smallskip
\noindent {\bf Convention}. For the curvature tensor of a connection
$\nabla$, we adopt the following definition
$R(X,Y)=\nabla_{[X,Y]}-[\nabla_X,\nabla_Y]$.

\smallskip

A straightforward computation gives the following relation between
$R^D$ and $R^g$ for $g\in {\mathfrak c}$. As in \cite{OT}, it is
convenient to set
\begin{equation}\label{Phi-g}
\Phi_g(X,Y)=(\nabla^g_{X}\varphi_g)(Y)+\frac{1}{2}\varphi_g(X)\varphi_g(Y)-\frac{1}{4}g(\varphi_g^{\sharp},\varphi_g^{\sharp})g(X,Y).
\end{equation}
Define an endomorphism $\Phi_g$ of $TM$ by
$g(\Phi_gX,Y)=\Phi_g(X,Y)$. Then
\begin{equation}\label{RD-Rg}
\begin{array}{r}
R^D(X,Y)Z=R^g(X,Y)Z+\displaystyle{\frac{1}{2}}\{\Phi_g(X,Y)Z-\Phi_g(Y,X)Z\\[8pt]
+\Phi_g(X,Z)Y-\Phi_g(Y,Z)X+g(X,Z)\Phi_gY-g(Y,Z)\Phi_gX\}.
\end{array}
\end{equation}
This formula implies the identities
\begin{equation}\label{Z-T}
g(R^D(X,Y)Z,T)+g(R^D(X,Y)T,Z)=d\varphi_g(X,Y)g(Z,T).
\end{equation}
\begin{equation}\label{XY-ZT}
\begin{array}{c}
2g(R^D(X,Y)Z,T)-2g(R^D(Z,T)X,Y)=d\varphi_g(X,Y)g(Z,T)
-d\varphi_g(Z,T)g(X,Y)\\[6pt]+d\varphi_g(X,Z)g(Y,T)+d\varphi_g(Y,T)g(X,Z)%\\[6pt]
-d\varphi_g(Y,Z)g(X,T)-d\varphi_g(X,T)g(Y,Z).
\end{array}
\end{equation}
\begin{equation}\label{Bianci}
R^D(X,Y)Z+R^D(Y,Z)X+R^D(Z,X)Y=0.
\end{equation}

\smallskip

\noindent {\bf Notation}. The Ricci tensor of the Weyl structure $({\mathfrak c},D)$ is defined by
$$
\rho_D(X,Z)=Trace_{g}\{Y\to g(R^D(X,Y)Z,Y)\},
$$
where $g\in{\mathfrak c}$.

If $J$ is a ${\mathfrak c}$-compatible almost-complex structure on
$M$, the $\ast$-Ricci tensor of the almost-Hermitian Weyl structure
$({\mathfrak c},D,J)$ is defined by
$$
\rho_{D}^{\ast}(X,Z)=Trace_{g}\{Y\to g(R^D(JY,X)JZ,Y)\}.
$$
Clearly, $\rho_{D}$ and $\rho_{D}^{\ast}$ do not depend on the particular choice of the metric $g$.

\smallskip

Identity (\ref{RD-Rg}) implies the following formulas.

\begin{prop}\label{Ric,star-Ric}
Let $\rho_g$ and $\rho_{g}^{\ast}$ be the Ricci tensor and the $\star$-Ricci tensor with respect to a Riemannian metric $g\in{\mathfrak c}$. Then
$$
\begin{array}{c}
\rho_{D}(X,Z)=\rho_{g}(X,Z)+\displaystyle{\frac{n-1}{2}(\nabla^g_{X}\varphi_g)(Z)-\frac{1}{2}(\nabla^g_{Z}\varphi_g)(X)}\\[10pt]
\hspace{2.6cm}-\displaystyle{\frac{n-2}{4}\big[||\varphi_g||^2g(X,Z)-\varphi_g(X)\varphi_g(Z)\big]-\frac{1}{2}(\delta^g\varphi_g)
g(X,Z)}.
\end{array}
$$
$$
\begin{array}{c}
\rho_{D}^{\ast}(X,Z)=\rho_{g}^{\ast}(X,Z)+(\nabla^g_{X}\varphi_g)(Z)-\displaystyle{\frac{1}{2}}\big[(\nabla^g_{Z}\varphi_g)(X)-(\nabla^g_{JX}\varphi_g)(JZ)\big]\\[10pt]
+\displaystyle{\frac{1}{4}}\big[\varphi_g(X)\varphi_g(Z)+\varphi_g(JX)\varphi_g(JZ)-||\varphi_g||^2g(X,Z)\big]\\[10pt]
+\displaystyle{\frac{1}{2}\big[\delta^g(J^{\ast}\varphi_g)-\varphi_g(\delta^gJ)\big]g(X,JZ)},
\end{array}
$$
where the norm and the codifferential are taken with respect to the
metric $g$.
\end{prop}
As is well-known, $\rho_{g}(X,Z)=\rho_{g}(Z,X)$ and
$\rho_{g}^{\ast}(X,Z)=\rho_{g}^{\ast}(JZ,JX)$.
Proposition~\ref{Ric,star-Ric} and these identities imply the
following.
\begin{cor}
$\rho_D(X,Z)-\rho_D(Z,X)=\displaystyle{\frac{n}{2}}d\varphi_{g}(X,Z)$.
\smallskip

$\rho_{D}^{\ast}(X,Z)-\rho_{D}^{\ast}(JZ,JX)=d\varphi_{g}(X,Z)+d\varphi_{g}(JX,JZ)$.
\end{cor}

\begin{cor}
$(1)$~The Ricci tensor $\rho_{D}$ is symmetric if and only if
$d\varphi_{g}=0$.

\smallskip

\noindent $(2)$ The $\ast$-Ricci tensor $\rho_{D}^{\ast}$ satisfies
the identity $\rho_{D}^{\ast}(X,Z)=\rho_{D}^{\ast}(JZ,JX)$ for all
$X,Z\in TM$ if and only if the $(1,1)$-component of $d\varphi_{g}$
w.r.t. $J$ vanishes.
\end{cor}
Claim $(1)$ is well-known. Note also that, by (\ref{change}), the
$2$-form $d\varphi_{g}$ does not depend on the choice of
$g\in{\mathfrak c}$.

\smallskip

Further, given a Riemannian metric $g$ on $M$, we shall often make use of the isomorphism $A(TM)\cong\Lambda^{2}TM$ that assigns to each $a\in A(T_{p}M)$ the 2-vector $a^{\wedge}$  determined by the identity
\begin{equation}\label{isom}
2g(a^{\wedge},X\wedge Y)=g(aX,Y),\quad X,Y\in T_{p}M,
\end{equation}
the metric on $\Lambda^{2}TM$ being defined by
$$
g(X_1\wedge X_2,X_3\wedge X_4)=\frac{1}{2}[g(X_1,X_3)g(X_2,X_4)-g(X_1,X_4)g(X_2,X_3)].
$$
Let $\widetilde g=e^fg$, $f$ being a smooth function, and denote by
$\widetilde{a^{\wedge}}$ the $2$-vector for which $\widetilde
g(\widetilde{a^{\wedge}},X\wedge Y) =\widetilde g(aX,Y)$. Then
$\widetilde{a^{\wedge}}=e^{-f}a^{\wedge}$.

\smallskip

Identity \cite[\rm{formula}~(11)]{D05} suggests the following.

\begin{lemma}\label{R[a,b]}
Let $g\in {\mathfrak c}$. Then, for every $a,b\in A(T_pM)$ and $X,Y\in T_pM$,
\begin{equation}\label{Rw}
\begin{array}{c}
G(R^D(X,Y)a,b)=g(R^D([a,b]^{\wedge})X,Y)\\[6pt]
-\displaystyle{\frac{1}{2}}[d\varphi_g([a,b]^{\wedge})g(X,Y)+d\varphi_g([a,b]X,Y)+d\varphi_g(X,[a,b]Y)].
\end{array}
\end{equation}
\end{lemma}

\begin{proof}

Let $E_1,...,E_n$ be a $g$-orthonormal basis of $T_pM$. Then
$$
[a,b]^{\wedge}=\frac{1}{2}\textstyle{\sum\limits_{i,j=1}^n}g([a,b]E_i,E_j)E_i\wedge E_j.
$$
Therefore, by (\ref{XY-ZT}),
$$
\begin{array}{l}
g(R^D([a,b]^{\wedge})X,Y)=\displaystyle{\frac{1}{2}}\sum\limits_{i,j=1}^ng(R^D(E_i,E_j)X,Y)g([a,b]E_i,E_j)\\[6pt]
=\displaystyle{\frac{1}{2}}\sum\limits_{i,j=1}^ng(R^D(X,Y)E_i,E_j)g([a,b]E_i,E_j)\\[8pt]
+\displaystyle{\frac{1}{2}}[d\varphi_g([a,b]^{\wedge})g(X,Y)+d\varphi_g([a,b]X,Y)+d\varphi_g(X,[a,b]Y)].
\end{array}
$$
Moreover,
$$
\begin{array}{l}
\sum\limits_{i,j=1}^ng(R^D(X,Y)E_i,E_j)g([a,b]E_i,E_j)\\[8pt]
=\sum\limits_{i,j=1}^n g(R^D(X,Y)E_i,E_j)[g(abE_i,E_j)+g(aE_i,bE_j)],\\[8pt]
=\sum\limits_{i=1}^ng(R^D(X,Y)E_i,abE_i)\\[6pt]
\hfill+\sum\limits_{i,j,k=1}^ng(R^D(X,Y)E_i,E_j)g(E_i,aE_k)g(E_j,bE_k)\\[8pt]
=-\sum\limits_{i=1}^n g(a(R^D(X,Y)E_i),bE_i)+\sum\limits_{k=1}^ng(R^D(X,Y)aE_k,bE_k),\\[8pt]
=2G(R^D(X,Y)a,b).
\end{array}
$$
This proves the lemma.
\end{proof}

\section{The induced connection $D'$ on the twistor space by a Weyl connection on the base manifold}

Let $I\in{\mathcal Z}$, and let ${\mathcal H}_I$ be the horizontal
subspace of $T_I A(TM)$ with respect to a Weyl connection $D$ on the
bundle $A(TM)$. Take a basis $e_1,...,e_{2m}$ of $T_pM$, $p=\pi(I)$,
such that  $Ie_{2k-1}=e_{2k}$, $Ie_{2k}=-e_{2k-1}$, $k=1,...,m$, and
$g(e_k,e_l)=\lambda_g\delta_{kl}$ for every $g\in {\mathfrak c}$,
where $\lambda_g$ is a positive constant depending on $g$. We can
extend this basis to a frame of vector fields $E_1,...,E_{2m}$ in a
(geodesically convex) neighbourhood of the point $p$ such that
$(E_i)_p=e_i$,  $D E_i|_p=0$ and $g(E_i,E_j)=f_g\delta_{ij}$,
$i,j=1,...,2m$, for every $g\in {\mathfrak c}$, where $f_g$ is a
positive smooth function. Define a section $K$ of $A(TM)$ by
$KE_{2k-1}=E_{2k}$, $KE_{2k}=-E_{2k-1}$. Obviously, $K$ takes values
in ${\mathcal Z}$, $K(p)=I$, and $D K|_p=0$. Hence the horizontal
space ${\mathcal H}_I=K_{\ast}(T_pM)$ is tangent to ${\mathcal Z}$.
Thus, we have the decomposition $T{\mathcal Z}={\mathcal
H}\oplus{\mathcal V}$ of the tangent bundle of ${\mathcal Z}$, where
${\mathcal V}=Ker (\pi_{\ast})$ is the vertical subbundle of
$T{\mathcal Z}$. The vertical space ${\mathcal V}_I$ at a point
$I\in {\mathcal Z}$ is the tangent space at $I$ of the fibre of
${\mathcal Z}$ through $I$. This fibre is the manifold of
${\mathfrak c}$-compatible complex structures on the vector space
$T_{\pi(I)}M$. Thus, the tangent space ${\mathcal V}_I$ consists of
skew-symmetric endomorphisms of $T_{\pi(I)}M$ anti-commuting with
$I$. In particular, if $V\in {\mathcal V}_I$ and $g\in {\mathfrak
c}$, then $g(IX,VX)=0$ for every $X\in T_{\pi(I)}M$. This fact and
identity (\ref{Z-T}) imply the following.
\begin{lemma}\label {J-V}
If $I\in{\mathcal Z}$ and $V\in{\mathcal V}_I$,
$$
G(R^D(X,Y)I,V)=-G(R^D(X,Y)V,I).
$$
\end{lemma}

Let $({\mathscr U},x_1,...,x_n)$ be a local coordinate system of $M$, and let $E_1,...,E_n$ be a frame of $TM$ on ${\mathscr U}$ such that,
for every metric $g\in {\mathfrak c}$, $g(E_i,E_j)=f_g\delta_{ij}$, $i,j=1,...,n$, where $f_g$ is a positive smooth function. Define sections $S_{ij}$ of $A(TM)$ by the formula
$$
S_{ij}E_l=\delta_{il}E_j - \delta_{lj}E_i,\quad i,j,l=1,...,n.
$$
These sections do not depend on the choice of the metric
$g\in{\mathfrak c}$ and $S_{ij}, i<j,$ form a $G$-orthonormal frame
of $A(TM)$. Set
$$
\tilde x_{i}(a)=x_{i}\circ\pi(a),\quad y_{jl}(a)=G(a,S_{jl}),\, j<l,
$$
for $a\in A(TM)$. Then $(\tilde x_{i},y_{jl})$ is a local coordinate system of the manifold $A(TM)$. For each vector field
$$
X=\sum_{i=1}^n X^{i}\frac{\partial}{\partial x_i}
$$
on ${\mathscr U}$, the horizontal lift $X^h$ on $\pi^{-1}({\mathscr U})$ is given by
\begin{equation}\label{Xh}
X^{h}=\sum_{i=1}^n (X^{i}\circ\pi)\frac{\partial}{\partial\tilde
x_i}-\sum_{j<l}\sum_{r<s}
y_{rs}(G(D_{X}S_{rs},S_{jl})\circ\pi)\frac{\partial}{\partial
y_{jl}}.
\end{equation}

Let $a\in A(TM)$ and $p=\pi(a)$. Then (\ref{Xh}) implies that, under the standard identification $T_{a}A(T_pM)\cong A(T_{p}M)$ (=the skew-symmetric endomorphisms of $T_{p}M$), we have the well-known formula
\begin{equation}\label{Lie-bra-hh}
[X^{h},Y^{h}]_{a}=[X,Y]^h_a + R^D(X,Y)a.
\end{equation}

Any (local) section $a$ of the bundle $A(TM)$ determines a (local)
vertical vector field $\widetilde a$ on ${\mathcal Z}$ defined by
$$
{\widetilde a}_I=\frac{1}{2}(a(q)+I\circ a(q)\circ I),\quad I\in
{\mathcal Z},\quad  q=\pi(I).
$$
Thus, if $aE_j=\sum\limits_{l=1}^n a_{jl}E_l$,
$$
\displaystyle{{\widetilde a}=\sum\limits_{j<l}\widetilde{a}_{jl}\frac{\partial}{\partial y_{jl}}},
$$
where
$$
\widetilde{a}_{jl}=\frac{1}{2}[a_{jl}\circ\pi+\textstyle{\sum\limits_{r,s=1}^n}y_{jr}(a_{rs}\circ\pi)y_{sl}].
$$
This formula and (\ref{Xh}) imply the following formula, which is a kind of folklore appearing in different contexts (see, for example, \cite[Lemma 1]{D05}, \cite[Appendice A]{Ga91}).

\begin{lemma}\label{Xh-a til}
If $I\in{\mathcal Z}$ and $X$ is a vector field on a neighbourhood of the point $p=\pi(I)$, then
$$
[X^h,\widetilde a]_I=(\widetilde{D_{X}a})_I.
$$
\end{lemma}

\begin{proof}
Take a frame $E_1,...,E_n$ in a neighbourhood of the point
$p=\pi(I)$ such that, for every metric $g\in {\mathfrak c}$,
$g(E_i,E_j)=f_g\delta_{ij}$, $i,j=1,...,n$, where $f_g$ is a
positive smooth function, and $DE_i|_p=0$ for $i=1,...,n$. Then
$(D_{X}a)(E_i)=\sum\limits_{l=1}^nX_p(a_{il})(E_l)_p$ and the result
follows by a trivial computation of $[X^h,\widetilde a]_I$ in the
coordinates $(\tilde x_{i},y_{jl})$ introduced above.
\end{proof}

\smallskip

Using the decomposition $T{\mathcal Z}={\mathcal H}\oplus {\mathcal
V}$, we can define a $1$-parameter family of Riemannian metrics on
${\mathcal Z}$ for every $g\in {\mathfrak c}$ setting
\begin{equation}\label{gt}
\tilde{g}_t(X^h+V,Y^h+W)=g(X,Y)+tG(V,W), \quad t>0,
\end{equation}
where $X^h,Y^h$ are the horizontal lifts of $X,Y\in TM$ (with respect to $D$) and $V,W$ are vertical vectors. Then the projection map $\pi:({\mathcal Z},\tilde{g}_t)\to (M,g)$ is a Riemannian submersion.

\smallskip

\noindent {\bf Notation}. Fix a metric $g\in {\mathfrak c}$ and
denote $\nabla^g$, $\varphi_g$ and $\delta^g$ simply by $\nabla$,
$\varphi$ and $\delta$. Let $\widetilde D=\widetilde D_{g,t}$ be the
Levi-Civita connection of the metric $\tilde g_t$.

\smallskip

\begin{lemma}\label{D-vv}
Let $I\in {\mathcal Z}$ and let $a$,$b$ be sections of $A(TM)$ in a neighbourhood of the point $p=\pi(I)$ such that $a(p)\in {\mathcal V}_I$ and $b(p)\in{\mathcal V}_I$. Then
$$
(\widetilde{D}_{\widetilde a}{\widetilde b})_I=0.
$$
\end{lemma}

\begin{proof}
In the notation of the preceding lemma and its proof, if $X\in T_pM$, then
$$
X^h_IG(\widetilde a,\widetilde b)=\sum\limits_{j,l=1}^n
X^h_I(\widetilde{a}_{jl}\widetilde{b}_{jl})=G((\widetilde{D_{X}a})_I,{\widetilde
b}_I)+G({\widetilde a}_I, (\widetilde{D_{X}b})_I).
$$
It follows from the latter identity, Lemma~\ref{Xh-a til} and the
Koszul formula for the Levi-Civita connection that $({\widetilde
D}_{\widetilde a}{\widetilde b})_I$ is orthogonal to every $X^h_I$,
so it is a vertical vector. Next, the Lie bracket of $\widetilde a$
and $\widetilde b$ is
$$
[\widetilde a,\widetilde b]_I=\frac{1}{4}\{[a(p),b(p)]I-I[a(p),b(p)]\},\quad p=\pi(I),
$$
where $[a(p),b(p)]$ is the commutator of the endomorphisms $a(p)$ and $b(p)$ of $T_pM$ and, as usual, juxtaposition means composition of maps. Since the endomorphisms $a(p)$ and $b(p)$ anti-commute with $I$,
$$
[a(p),b(p)]I-I[a(p),b(p)]=0.
$$
Let $W\in {\mathcal V}_I$ and let $c$ be a local section of $A(TM)$ such that $c(p)=W$. Then
$$
\begin{array}{c}
{\widetilde c}_I(G(\widetilde a,\widetilde b))=\sum\limits_{j,l=1}^n W(\widetilde{a}_{jl}\widetilde{b}_{jl})\\[6pt]
=\frac{1}{4}\{G([W,a(p)]I-I[W,a(p)],b(p))+G([W,b(p)]I-I[W,b(p)],a(p))\}.
\end{array}
$$
Hence
$$
{\widetilde c}_I(G(\widetilde a,\widetilde b))=0.
$$
Now, it follows from the Koszul formula that
$\tilde{g}_t((\widetilde{D}_{\widetilde a}{\widetilde
b})_I,{\widetilde c}_I)=0$. This proves the lemma.
\end{proof}

Let $\{V_{\alpha}\}$, $\alpha=1,...,m^2-m$, be a basis of a vertical space ${\mathcal V}_I$. Take local sections $a_{\alpha}$ of $A(TM)$ such that $a_{\alpha}(p)=V_{\alpha}$, $p=\pi(I)$.
Then the vertical vector fields $\{ \widetilde{a_{\alpha}}\}$ constitute a frame in a neighbourhood of $I$.
It follows from Lemma~\ref{D-vv} that the covariant derivative of a vertical vector field in a vertical direction is a vertical vector.
Thus, the fibres of the bundle $\pi:({\mathcal Z},\widetilde g_t)\to (M,g)$ are totally geodesic submanifolds. This, of course, follows also from  the Vilms theorem
(see, for example, \cite[Theorem 9.59]{Besse}).

\smallskip

The Koszul formula, identities (\ref{Lie-bra-hh}) and (\ref{Weyl}), and the fact that the fibres of the twistor bundle are totally geodesic submanifolds imply the following formulas.

\begin{lemma}\label{LC}
If $X,Y$ are vector fields on $M$ and $V$ is a vertical vector field on ${\mathcal Z}$, then at $I\in{\mathcal Z}$
\begin{equation}\label{D-hh}
(\widetilde
D_{X^h}Y^h)_{I}=(D_{X}Y)^h_{I}+\frac{1}{2}[\varphi(X)Y^h_{I}+\varphi(Y)X^h_{I}-g(X,Y)(\varphi^{\sharp})^h_{I}]+\frac{1}{2}R^D(X,Y)I.
\end{equation}
Also, $\widetilde D_{V}X^h={\mathcal H} \widetilde D_{X^h}V$, where
${\mathcal H}$ means the horizontal component. Moreover,
\begin{equation}\label{D-vh}
\widetilde g_t(\widetilde
D_{V}X^h,Y^h)_{I}=-\frac{t}{2}G(R^D_{\pi(I)}(X,Y)I,V).
\end{equation}
\end{lemma}

\smallskip

\noindent {\bf Notation}. The lemma above suggests defining a new
torsion-free connection $D'={D'}_{g,t}$ on ${\mathcal Z}$ as
follows. Let $I\in {\mathcal Z}$. Let $X,Y$ be vector fields on $M$
in a neighbourhood of $p=\pi(I)$, and let $V,W$ be vertical vector
fields on ${\mathcal Z}$ in a neighbourhood of $I$. Set
\begin{equation}\label{new conn}
\begin{array}{c}
(D'_{X^h}Y^h)_I=(D_{X}Y)^h_I+\frac{1}{2}R^D(X,Y)I, \\[6pt]
D'_{V}X^h=\widetilde D_{V}X^h, \quad D'_{X^h}V=\widetilde D_{X^h}V,\quad D'_{V}W=\widetilde D_{V}W.
\end{array}
\end{equation}

Recall that $R^D(X,Y)I$ in the first formula is a vertical vector at $I$.
Clearly, the fibres of the bundle $\pi:{\mathcal Z}\to M$ are still totally geodesic submanifolds when ${\mathcal Z}$ is considered with the connection $D'$.

\section{The second fundamental form w.r.t. $D$ and $D'$ of an almost-Hermitian structure as a map into the twistor space}

Let $J$ be a ${\mathfrak c}$-compatible almost-complex structure. We
endow $M$ with the orientation yielded by $J$. Then $J$ can be
considered as a section of the positive twistor space of
$(M,{\mathfrak c})$.

\smallskip

\noindent {\bf Notation}. Henceforth, $\pi:{\mathcal Z}\to M$ will denote the positive twistor space of $(M,{\mathfrak c})$.

\smallskip

Let $J^{\ast}T{\mathcal Z}\to M$ be the pull-back of the bundle $T{\mathcal Z}\to {\mathcal Z}$ under the map $J:M\to{\mathcal Z}$. Then we can consider the differential $J_{\ast}:TM\to T{\mathcal Z}$ as a section of the bundle $Hom(TM,J^{\ast}T{\mathcal Z})\to M$.

\smallskip

\noindent {\bf Notation}. Denote by $D^{\ast}$ the connection on
$J^{\ast}T{\mathcal Z}$ induced by the torsion-free connection $D'$
on $T{\mathcal Z}$ defined in the preceding section. The Weyl
connection $D$ on $TM$ and the connection $D^{\ast}$ on
$J^{\ast}T{\mathcal Z}$ induce a connection $\widehat{D}$ on the
bundle $Hom(TM,J^{\ast}T{\mathcal Z})$. By definition, the second
fundamental form of the map $J:M\to{\mathcal Z}$ with respect to the
connections $D$ and $D'$ is the bilinear form defined as
$$
II_{J}(D,D')(X,Y)=(\widehat D_{X}J_{\ast})(Y),\quad X,Y\in TM.
$$
This form is symmetric since $D$ and $D'$ are torsion-free.
Following \cite{K09}, we say that the map $J:M\to{\mathcal Z}$ is
$(D,D')$-pseudo harmonic if
$$
Trace_{g} II_{J}(D,D')=0,
$$
where $g\in{\mathfrak c}$ and the trace clearly does not depend on the choice of $g$.

To compute this trace, we need the following lemma.

\begin{lemma}\label{v-frame}
For every $p\in M$, there exists a $\widetilde g_t$-orthonormal
frame of vertical vector fields
$\{V_{\alpha}:~\alpha=1,...,m^2-m\}$, $2m=dim\,M$, in a
neighbourhood of the point $J(p)$ such that

\noindent $(1)$ $\quad (D'_{V_{\alpha}}V_{\beta})_{J(p)}=0$,~\quad~$\alpha,\beta=1,...,m^2-m$.

\noindent $(2)$ $\quad$ If $X$ is a vector field near the point $p$, $[X^h,V_{\alpha}]_{J(p)}=0$.

\noindent $(3)$ $\quad$  $D_{X_p}(V_{\alpha}\circ J)\perp {\cal V}_{J(p)}$, where $V_{\alpha}\circ J$ is considered as a section of $A(TM)$.
\end{lemma}

The proof goes along the same lines as the proof of \cite[Lemma 8]{D17}, so it will be skipped.

\begin{prop}\label{covder-dif}
For every  $X,Y,Z\in T_pM$, $p\in M$, and $W\in{\mathcal V}_{J(p)}$

$$
\begin{array}{c}
2\widetilde g_t((\widehat D_{X} J_{\ast})(Y),Z^h_{J(p)})
=-tG(R^D(Y,Z)J(p),D_{X}J)-tG(R^D(X,Z)J(p),D_{Y}J),
\end{array}
$$
$$
\begin{array}{c}
2\widetilde g_t((\widehat D_{X}J_{\ast})(Y),W) =tG(D^{2}_{XY}J+D^{2}_{YX}J,W),
\end{array}
$$
where $D^{2}_{XY}J=D_XD_Y J-D_{D_XY}J$ is the second covariant derivative of $J$ considered as a section of $A(TM)$.
\end{prop}

\begin{proof} Extend $X$ and $Y$ to vector fields in a neighbourhood of the point $p$. Let $V_1,...,V_{m^2-m}$ be a $\widetilde g_t$-orthonormal frame of vertical vector fields with the properties $(1)$ - $(3)$ stated in Lemma~\ref{v-frame}.

We have
$$
J_{\ast}\circ Y=Y^h\circ J+D_{Y}J=Y^h\circ J+\textstyle{\sum\limits_{\alpha=1}^{m^2-m}}\widetilde g_t(D_{Y}J,V_\alpha\circ J)(V_{\alpha}\circ J),
$$
hence
$$
\begin{array}{c}
D^{\ast}_{X}(J_{\ast}\circ Y)=(D'_{J_{\ast}X}Y^h)\circ J+
\sum\limits_{\alpha=1}^{m^2-m}\widetilde g_t(D_{Y}J,V_{\alpha}\circ J)(D'_{J_{\ast}X}V_{\alpha})\circ J\\[8pt]
+t\sum\limits_{\alpha=1}^{m^2-m}G(D_{X}D_{Y} J,V_{\alpha}\circ J)(V_{\alpha}\circ J).
\end{array}
$$
By (\ref{D-hh}) and (\ref{new conn}),
$$
\begin{array}{c}
(D'_{J_{\ast}X}Y^h)_{J(p)}= (D_{X}Y)^h_{J(p)}
%+\displaystyle{\frac{1}{2}}[\varphi(X)Y^h_{J(p)}+\varphi(Y)X^h_{J(p)}-g(X,Y)(\varphi^{\sharp})^h_{J(p)}]\\[8pt]
+\displaystyle{\frac{1}{2}}R^D(X\wedge Y)J(p)+(D'_{D_{X}J}Y^h)_{J(p)}.
\end{array}
$$
By Lemma~\ref{v-frame},
$$
\begin{array}{c}
(D'_{J_{\ast}X}V_{\alpha})_{J(p)}=D'_{X^h_{J(p)}}V_{\alpha}+D'_{D_{X_p}J}V_{\alpha},\\[8pt]
=[X^h,V_{\alpha}]_{J(p)}+(D'_{V_{\alpha}}X^h)_{J(p)}+D'_{D_{X_p}J}V_{\alpha}=(D'_{V_{\alpha}}X^h)_{J(p)}.
\end{array}
$$
Thus,
$$
\begin{array}{c}
{D}^{\ast}_{X_p}(J_{\ast}\circ Y)= (D_{X}Y)^h_{J(p)}
%+\displaystyle{\frac{1}{2}}[\varphi(X)Y^h_{J(p)}+\varphi(Y)X^h_{J(p)}-g(X,Y)(\varphi^{\sharp})^h_{J(p)}]\\[8pt]
+(D'_{D_{X}J}Y^h)_{J(p)}+(D'_{D_{Y}J}X^h)_{J(p)}\\[8pt]
+\displaystyle{\frac{1}{2}}R^D(X\wedge Y)J(p)+t\sum\limits_{\alpha=1}^{m^2-m}G(D_{X_p}D_{Y}J,V_{\alpha}\circ J)_{p}(V_{\alpha})_{J(p)},
\end{array}
$$
where
$$
R^D(X\wedge Y)J(p)=D_{D_{X_p}Y}J-D_{D_{Y_p}X}J-D_{X_p}D_{Y}J+D_{Y_p}D_{X}J
$$
by the definition we have adopted, and
$$
R^D(X\wedge Y)J(p)=t\sum\limits_{\alpha=1}^{m^2-m}G(R^D(X,Y)J(p),(V_{\alpha})_{J(p)})(V_{\alpha})_{J(p)}.
$$
It follows that
$$
\begin{array}{c}
(\widehat D_{X}J_{\ast})(Y)=D^{\ast}_{X_p}(J_{\ast}\circ Y)-J_{\ast}(D_{X}Y),\\[8pt]
=D^{\ast}_{X_p}(J_{\ast}\circ Y)-(D_{X}Y)^h_{J(p)}-D_{D_{X_p}Y}J,\\[8pt]
=(D'_{D_{X}J}Y^h)_{J(p)}+(D'_{D_{Y}J}X^h)_{J(p)}\\[8pt]
+\displaystyle{\frac{1}{2}}{\cal V}(D_{X_p}D_{Y}J-D_{D_{X_p}Y}J+D_{Y_p}D_{X}J-D_{D_{Y_p}X}J).
\end{array}
$$

Now, the proposition follows from (\ref{D-vh}).
\end{proof}

\begin{lemma}\label{R-[J,DJ]}
For every $X,Y,Z\in T_pM$,
$$
\begin{array}{c}
G(R^D(X,Z)J,D_{Y}J)=2g(R^D((J\nabla_{Y}J)^{\wedge})X,Z)-g(R^D(\varphi^{\sharp}\wedge
Y-J\varphi^{\sharp}\wedge JY)X,Z)\\[8pt]
-d\varphi((J\nabla_{Y}J)^{\wedge})g(X,Z)-d\varphi((J\nabla_{Y}J)(X),Z)-d\varphi(X,(J\nabla_{Y}J)(Z))\\[6pt]
+\displaystyle{\frac{1}{2}}d\varphi(\varphi^{\sharp}\wedge
Y-J\varphi^{\sharp}\wedge
JY)g(X,Z)\\[8pt]
+\displaystyle{\frac{1}{2}}[\varphi(JX)d\varphi(JY,Z)+\varphi(X)d\varphi(Y,Z)-g(Y,JX)d\varphi(J\varphi^{\sharp},Z)-g(Y,X)d\varphi(\varphi^{\sharp},Z)]\\[8pt]
+\displaystyle{\frac{1}{2}}[\varphi(JZ)d\varphi(X,JY)+\varphi(Z)d\varphi(X,Y)-g(Y,JZ)d\varphi(X,J\varphi^{\sharp})-g(Y,Z)d\varphi(X,\varphi^{\sharp})].
\end{array}
$$
\end{lemma}

\begin{proof} By Lemma~\ref{R[a,b]},
$$
\begin{array}{c}
G(R^D(X,Z)J,D_{Y}J)=g(R^D([J,D_{Y}J]^{\wedge})X,Z)\\[6pt]
-\displaystyle{\frac{1}{2}}[d\varphi([J,D_{Y}J]^{\wedge})g(X,Z)+d\varphi([J,D_{Y}J]X,Z)+d\varphi(X,[J,D_{Y}J]Z)],\\[6pt]

=2g(R^D((JD_{Y}J)^{\wedge})X,Z) \\[6pt]
\hspace{1.4cm}-d\varphi((JD_{Y}J)^{\wedge})g(X,Z)-d\varphi((JD_{Y}J)X,Z)-d\varphi(X,(JD_{Y}J)Z).
\end{array}
$$

By (\ref{Weyl}),
\begin{equation}\label{DJ}
J(D_{Y}J)(E)=J(\nabla_{Y}J)(E)-\frac{1}{2}[\varphi(JE)JY+\varphi(E)Y-g(Y,JE)J\varphi^{\sharp}-g(Y,E)\varphi^{\sharp}],
\end{equation}

for every $Y,E\in T_pM$. It follows that

$$
(JD_{Y}J)^{\wedge}=(J\nabla_{Y}J)^{\wedge}-\frac{1}{2}(\varphi^{\sharp}\wedge Y-J\varphi^{\sharp}\wedge JY).
$$
Now, the lemma follows easily.
\end{proof}

Proposition~\ref{covder-dif}, Lemma~\ref{R-[J,DJ]} and identity (\ref{Z-T}) imply the following.
\begin{cor}\label{H-tr}
For every tangent vector $Z\in T_pM$,
$$
\begin{array}{c}
-\frac{1}{t}\widetilde g_t({\mathcal H}(Trace_g \widehat D J_{\ast}),Z^h_{J(p)})\\[6pt]

=2Trace_g\{T_pM\ni X\to g(R^D((J\nabla_{X}J)^{\wedge})X,Z)\}+\rho_D(\varphi^{\sharp},Z)-\rho_D^{\ast}(J\varphi^{\sharp},JZ)\\[6pt]

-d\varphi((J\nabla_{Z}J)^{\wedge})+d\varphi(J\delta J,Z)-Trace_g\{T_pM\ni X\to d\varphi(X,(J\nabla_{X}J)(Z))\} \\[6pt]

+\varphi(JZ)d\varphi (J^{\wedge})\\[6pt]

-\displaystyle{(\frac{n}{2}-1)}d\varphi(\varphi^{\sharp},Z)+d\varphi(J\varphi^{\sharp},JZ).
\end{array}
$$
\end{cor}

\begin{lemma}
For every $Z,U\in T_pM$,
$$
\begin{array}{c}
g(Trace\{X\to (D^2_{XX}J)(Z)\},U)=g(Trace\{X\to (\nabla^2_{XX}J)(Z)\},U)\\[10pt]

+\displaystyle{\frac{1}{2}(2-n)g((\nabla_{\varphi^{\sharp}}J)(Z),U)-\varphi(Z)g(\delta J,U)+\varphi(U)g(\delta J,Z)}\\[10pt]

+\displaystyle{g((\nabla_{Z}J)(U),\varphi^{\sharp})-g((\nabla_{U}J)(Z),\varphi^{\sharp})+\frac{1}{2}d\varphi(JZ\wedge U+Z\wedge JU)}\\[10pt]

+\displaystyle{\frac{1}{4}(n-3)[\varphi(Z)\varphi(JU)-\varphi(JZ)\varphi(U)]+\frac{1}{2}||\varphi^{\sharp}||^2g(Z,JU)}.
\end{array}
$$
\end{lemma}

\begin{proof} We have
$$
\begin{array}{l}
g((D^2_{XY}J)(Z),U)=g((D_{X}(D_{Y}J))(Z),U)-g((D_{Y}J)(D_{X}Z),U)\\[6pt]
\hfill -g((D_{D_{X}Y}J)(Z),U),\\[10pt]

=X.g((D_{Y}J)(Z),U)-g((D_{Y}J)(Z),D_{X}U)-\varphi(X)g((D_{Y}J)(Z),U)\\[6pt]
\hfill -g((D_{Y}J)(D_{X}Z),U) -g((D_{D_{X}Y}J)(Z),U).\\[10pt]
\end{array}
$$

Also, by (\ref{Weyl}),
$$
g((D_{Y}J)(Z),U)=g((\nabla_{Y}J)(Z),U)-g(JY\wedge\varphi^{\sharp}+Y\wedge J\varphi^{\sharp},Z\wedge U).
$$
Applying these identities and identity (\ref{Weyl}), after a simple computation, we obtain
$$
\begin{array}{l}
g((D^2_{XY}J)(Z),U)=g((\nabla^2_{XY}J)(Z),U)\\[10pt]

+\displaystyle{\frac{1}{2}}[\varphi(X)g((\nabla_{Y}J)(Z),U)+\varphi(Y)g((\nabla_{X}J)(Z),U)\\[6pt]

\hfill+\varphi(Z)g((\nabla_{Y}J)(X),U)-\varphi(U)g((\nabla_{Y}J)(X),Z)]\\[10pt]

+\displaystyle{\frac{1}{2}}[g(X,U)g((\nabla_{Y}J)(\varphi^{\sharp}),Z)-g(X,Z)g((\nabla_{Y}J)(\varphi^{\sharp}),U)\\[6pt]
\hfill-g(X,Y)g((\nabla_{\varphi^{\sharp}}J)(Z),U)]\\[10pt]

-g((\nabla_{X}J)(Y)\wedge\varphi^{\sharp}+Y\wedge(\nabla_{X}J)(\varphi^{\sharp})+JY\wedge\nabla_{X}\varphi^{\sharp}\\[6pt]

\hfill+Y\wedge J\nabla_{X}\varphi^{\sharp},Z\wedge U)\\[6pt]

+\displaystyle{\frac{1}{2}}g(\varphi(JY)X\wedge\varphi^{\sharp}+||\varphi^{\sharp}||^2JY\wedge X+\varphi(X)Y\wedge J\varphi^{\sharp} -\varphi(Y)JX\wedge\varphi^{\sharp}\\[6pt]

\hfill-g(X,JY)Y\wedge\varphi^{\sharp}+g(X,Y)J\varphi^{\sharp}\wedge\varphi^{\sharp},Z\wedge U).
\end{array}
$$

Finally, note that

$$
\begin{array}{c}
g(\nabla_{JZ}\varphi^{\sharp},U) -g(\nabla_{JU}\varphi^{\sharp},Z)+g(\nabla_{Z}\varphi^{\sharp},JU)-g(\nabla_{U}\varphi^{\sharp},JZ)\\[6pt]

=(\nabla_{JZ}\varphi)(U)-(\nabla_{JU}\varphi)(Z)+(\nabla_{Z}\varphi)(JU)-(\nabla_{U}\varphi)(JZ),\\[6pt]

=d\varphi(JZ\wedge U)+d\varphi(Z\wedge JU).
\end{array}
$$
\end{proof}

Let $I\in {\mathcal Z}$. Take a metric $g_0\in{\mathfrak c}$. Let
$e_1,...,e_{n}$ be a $g_0$-orthonormal basis of $T_pM$, $p=\pi(I)$,
such that $Ie_{2k-1}=e_{2k}$, $k=1,...,m$. Extend this basis to a
frame of vector fields $E_1,...,E_n$ in a neighbourhood of $p$ such
that $(E_i)_p=e_i$, $DE_i|_p=0$, $i=1,...,n$, and, for every
$g\in{\mathfrak c}$, $g(E_k,E_l)=f_g\delta_{kl}$, $k,l=1,...,n$,
where $f_g$ is a smooth positive function. As above, define sections
$S_{ij}$, $1\leq i,j\leq n$, of $A(TM)$, by
$$
S_{ij}E_l=\delta_{il}E_j - \delta_{lj}E_i,\quad i,j,l=1,...,n.
$$
Set
$$
\begin{array}{c}
A_{r,s}=\frac{1}{\sqrt 2}(S_{2r-1,2s-1}-S_{2r,2s}),\quad B_{r,s}=\frac{1}{\sqrt 2}(S_{2r-1,2s}+S_{2r,2s-1}),\\[6pt]
r=1,...,m-1,\>s=r+1,...,m.
\end{array}
$$
The endomorphisms $\{(A_{r,s})_p,(B_{r,s})_p\}$ of $T_pM$,
$r=1,...,m-1$, $s=r+1,...,m$, constitute a $G$-orthonormal basis of
the vertical space ${\mathcal V}_I$. It follows that identifying
$A(T_pM)$ with $\Lambda^2T_pM$ by means of the isomorphism $a\to
a^{\wedge}$ given by (\ref{isom}),  every vector of ${\mathcal V}_I$
is a linear combination of vectors of the form $Z\wedge U-JZ\wedge
JU$, $Z,U\in T_pM$. Thus, $Trace\{X\to G(D^2_{XX}J,W)\}=0$ for every
vertical vector $W$ of ${\mathcal Z}$ if and only if $Trace\{X\to
g((D^2_{XX}J)^{\wedge},Z\wedge U-JZ\wedge JU)\}=0$ for every tangent
vectors $Z,U$ of $M$. In view of this remark and
Proposition~\ref{covder-dif}, we state the following.

\begin{cor}\label{V-tr}
If $X,Z,U\in T_pM$,
$$
\begin{array}{c}

g((Trace\,D^2_{XX}J)(Z),U)-g((Trace\,D^2_{XX}J)(JZ),JU)\\[10pt]

=g((Trace\,\nabla^2_{XX}J)(Z),U)-g((Trace\,\nabla^2_{XX}J)(JZ),JU)\\[10pt]

+(2-n)g((\nabla_{\varphi^{\sharp}}J)(Z),U)\\[10pt]

-\varphi(Z)g(\delta J,U)+\varphi(JZ)g(\delta J,JU)+\varphi(U)g(\delta J, Z)-\varphi(JU)g(\delta J, JZ)\\[10pt]

-g((\nabla_{Z}J)(\varphi^{\sharp}),U)+g((\nabla_{U}J)(\varphi^{\sharp}),Z)+g((\nabla_{JZ}J)(\varphi^{\sharp}),JU)\\[10pt]

-g((\nabla_{JU}J)(\varphi^{\sharp}),JZ)+d\varphi(JZ\wedge U+Z\wedge JU).

\end{array}
$$
\end{cor}

\section{Hermitian structures on Weyl manifolds yielding pseudo-harmonic maps into the twistor space}

\smallskip

\noindent {\bf Notation}. Denote by $N$ the Nijenhuis tensor of $J$:
$$
N(Y,Z)=-[Y,Z]+[JY,JZ]-J[Y,JZ]-J[JY,Z].
$$

Let $\Omega(X,Y)=g(JX,Y)$ be the fundamental $2$-form of the
almost-Hermitian manifold $(M,g,J)$ ($g\in{\mathfrak c}$ being a
fixed metric as in the preceding sections).

\smallskip

It is well-known (and easy to check) that
\begin{equation}\label{nJ}
2g((\nabla_XJ)(Y),Z)=d\Omega(X,Y,Z)-d\Omega(X,JY,JZ)+g(N(Y,Z),JX),
\end{equation}
for all $X,Y,Z\in TM$. The condition that $J$ is integrable ($N=0$) is equivalent to $(\nabla_{X}J)(Y)=(\nabla_{JX}J)(JY)$, $X,Y\in TM$
\cite[Corollary 4.2]{G65}.

Suppose that the almost-complex structure $J$ is integrable. Let
$$
\theta=-\frac{2}{n-2}\delta\Omega\circ J
$$
be the Lee form of the Hermitian structure $(g,J)$. If $n=4$, the
Lee form satisfies the identity $d\Omega=\theta\wedge\Omega$. For
$n\geq 6$, this identity is satisfied if and only $(g,J)$ is a
locally conformally K\"ahler structure. As is well-known, it implies
$d\theta=0$ when $n\geq 6$. Note also  that, in any dimension, an
Hermitian structure $(g,J)$ is  locally conformally K\"ahler
 if and only if $d\Omega=\theta\wedge\Omega$ and $d\theta=0$, see, for example,
 \cite{DO,V76}. Clearly, if $(g,J)$ is such a structure, then for any metric $g^1$
conformal to $g$, the structure $(g^1,J)$ is also locally
conformally K\"ahler.

\smallskip

\noindent {\it Assumption}. We assume that the almost-complex
structure $J$ is integrable and the Lee from satisfies the identity
\begin{equation}\label{dOm}
d\Omega=\theta\wedge\Omega.
\end{equation}

\smallskip

\noindent {\bf Remark}. If the latter identity is satisfied by an
almost-Hermitian structure $(g,J)$, it is also satisfied by any
almost-Hermitian structure $(g^1,J)$ with $g^1=e^fg$, $f$ being a
smooth function,  since $\Omega^1=e^f\Omega$ and
$\theta^1=\theta+df$ by the first identity of (\ref{change}).

\smallskip

Let $B$ be the vector field on $M$ dual to the Lee form $\theta=-\delta\Omega\circ J$ with respect to the metric $g$. Thus
$$
B=\frac{2}{n-2}J\delta J.
$$

Identities (\ref{nJ}) and (\ref{dOm}) imply the following well-known formula
\begin{equation}\label{nablaJ}
2(\nabla_XJ)(Y)=g(JX,Y)B-g(B,Y)JX+g(X,Y)JB-g(JB,Y)X.
\end{equation}

Now, we are going to apply Corollary~\ref{H-tr}.

\smallskip

Formula (\ref{nablaJ}) implies that, for $X,Y,Z\in T_pM$,
$$
g((J\nabla_{X}J)^{\wedge},Y\wedge Z)=\frac{1}{2}g(J(\nabla_{X}J)(Y),Z)=\frac{1}{2}g(B\wedge X-JB\wedge JX,Y\wedge Z).
$$
So,
$$
(J\nabla_{X}J)^{\wedge}=\frac{1}{2}(B\wedge X-JB\wedge JX),
$$
and applying (\ref{Z-T}), we get
$$
\begin{array}{c}
2Trace\{TM\ni X\to g(R^D((J\nabla_{X}J)^{\wedge})X,Z)\}\\[6pt]
=-\rho_{D}(B,Z)+\rho_{D}^{\ast}(JB,JZ)+d\varphi(B,Z)-d\varphi(JB,JZ).
\end{array}
$$
Moreover, by (\ref{nablaJ}),
$$
\begin{array}{c}
Trace\{TM\ni X\to d\varphi(X,(J\nabla_{X}J)(Z))\}\\[6pt]
=\displaystyle{\frac{1}{2}}[d\varphi(B,Z)+d\varphi(JB,JZ)]+\theta(JZ)d\varphi(J^{\wedge}).
\end{array}
$$

It follows from  Corollary~\ref{H-tr} that
$$
\begin{array}{c}

-\frac{1}{t}\widetilde g_t({\mathcal H}(Trace_g \widehat D J_{\ast}),Z^h_{J(p)})\\[6pt]

=\displaystyle{(\frac{n}{2}-1)}d\varphi((\theta-\varphi)^{\sharp},Z)-d\varphi(J(\theta-\varphi)^{\sharp},JZ)-(\theta-\varphi)(JZ)d\varphi(J^{\wedge})\\[6pt]
\hfill-\rho_{D}((\theta-\varphi)^{\sharp},Z)+\rho^{\ast}_{D}(J(\theta-\varphi)^{\sharp},JZ).

\end{array}
$$

Thus, we have the following statement.

\begin{lemma}\label{H-comp-int-Weyl}
Suppose that $J$ is a Hermitian structure whose Lee from satisfies
identity (\ref{dOm}). Then ${\mathcal H}(Trace\,\widehat D^2
J_{\ast})=0$ if and only if for every $Z\in TM$
$$
\begin{array}{c}
\displaystyle(\frac{n}{2}-1)d\varphi((\theta-\varphi)^{\sharp},Z)-d\varphi(J(\theta-\varphi)^{\sharp},JZ)-(\theta-\varphi)(JZ)d\varphi(J^{\wedge})\\[6pt]

\hfill-\rho_{D}((\theta-\varphi)^{\sharp},Z)+\rho^{\ast}_{D}(J(\theta-\varphi)^{\sharp},JZ)=0.
\end{array}
$$
\end{lemma}

Next, by Corollary~\ref{V-tr} and identity (\ref{nablaJ}),
$$
\begin{array}{c}
g((Trace\,D^2_{XX}J)(Z),U)-g((Trace\,D^2_{XX}J)(JZ),JU) \\[6pt]

=g((Trace\,\nabla^2_{XX}J)(Z),U)-g((Trace\,\nabla^2_{XX}J)(JZ),JU)\\[6pt]

+d\varphi(JZ\wedge U+Z\wedge JU).
\end{array}
$$

It follows from (\ref{nablaJ}) that

$$
\begin{array}{c}

2(\nabla^{2}_{XY}J)(Z)=2\big[\nabla_{X}(\nabla_{Y}J)(Z)-(\nabla_{Y}J)(\nabla_{X}Z)-(\nabla_{\nabla_{X}Y}J)(Z)\big]\\[6pt]

=g((\nabla_{X}J)(Y),Z)B-g(B,Z)(\nabla_{X}J)(Y)+g(Y,Z)(\nabla_{X}J)(B)\\[6pt]

-g((\nabla_{X}J)(B),Z)Y+g(JY,Z)\nabla_{X}B-g(\nabla_{X}B,Z)JY\\[6pt]

+g(Y,Z)J\nabla_{X}B-g(J\nabla_{X}B,Z)Y

\end{array}
$$

for every $X,Y,Z\in TM$. This identity and the identity (\ref{nablaJ}) imply

$$
\begin{array}{c}

2g((Trace\,D^2_{XX}J)(Z),U)=\displaystyle{\frac{n-4}{2}}[g(JB,Z)g(B,U)-g(JB,U)g(B,Z)]\\[6pt]

+ d\theta(U,JZ)+d\theta(JU,Z)-||B||^2g(JZ,U).

\end{array}
$$

We also have

$$
\begin{array}{c}

(2-n)g((\nabla_{\varphi^{\sharp}}J)(Z),U)\\[6pt] -\varphi(Z)g(\delta J,U)+\varphi(JZ)g(\delta J,JU)+\varphi(U)g(\delta J,Z)-\varphi(JU)g(\delta J, JZ) \\[6pt]

=\displaystyle{\big[\frac{2-n}{2}+\frac{2}{n-2}\big]}[-\varphi(JZ)\theta(U)+\varphi(JU)\theta(Z)

\hfill -\varphi(Z)\theta(JU)+\varphi(U)\theta(JZ)]\\[6pt]

=\displaystyle{\frac{n(n-4)}{2(n-2)}}(\varphi\wedge\theta )(JZ\wedge U+Z\wedge JU).
\end{array}
$$

Hence, by Corollary~\ref{V-tr},
$$
\begin{array}{c}

g((Trace\,D^2_{XX}J)(Z),U)-g((Trace\,D^2_{XX}J)(JZ),JU) \\[10pt]

=(d(\varphi-\theta)+\displaystyle{\frac{n(n-4)}{2(n-2)}}\varphi\wedge\theta)(JZ\wedge U+Z\wedge JU).

\end{array}
$$
As we have remarked, ${\mathcal V}(Trace\,\widehat D^2 J_{\ast})=0$
if and only if $g((Trace\,D^2J)(Z),U)-g((Trace\,D^2J)(JZ),JU)=0$ for
all $Z,U$. Note also that a $2$-form $\alpha$ is of type $(1,1)$
with respect to $J$ exactly when $\alpha(JZ\wedge U+Z\wedge JU)=0$
for every $Z,U\in TM$. Thus, we have the following.

\begin{lemma}\label{V-comp-int-Weyl}
Suppose that $J$ is a Hermitian structure whose Lee from satisfies
identity (\ref{dOm}). Then ${\mathcal V}(Trace\,\widehat D^2
J_{\ast})=0$ if and only if the $2$-form
$$
d(\varphi-\theta)+\frac{n(n-4)}{2(n-2)}\varphi\wedge\theta
$$
is of type $(1,1)$ with respect to $J$.
\end{lemma}

Proposition~\ref{covder-dif}, Lemmas~\ref{H-comp-int-Weyl} and \ref{V-comp-int-Weyl} imply the following.

\begin{theorem}\label{psd-harm-int}
Let $(M,{\mathfrak c})$ be a conformal $n$-dimensional  manifold
with a Weyl connection $D$,  let $g\in{\mathfrak c}$ and set
$\varphi=\varphi_{g}$. Suppose that $J$ is a Hermitian structure on
$(M,{\mathfrak c})$ such that Lee form $\theta$ of $(g,J)$ satisfies
identity (\ref{dOm}). Let $D'$ be the torsion-free connection on
${\mathcal Z}$ defined by (\ref{new conn}). Then $J$ determines a
$(D,D')$-pseudo-harmonic map from $M$ into the positive twistor
space ${\mathcal Z}$ of $(M,{\mathfrak c})$ if and only if the
following two conditions are satisfied.

\smallskip

(i) The $2$-form
$$
d(\theta-\varphi)+\frac{n(n-4)}{2(n-2)}\theta\wedge\varphi
$$
is of type $(1,1)$ with respect to $J$.

\smallskip

(ii) For every tangent vector $Z\in TM$,
$$
\begin{array}{c}

\displaystyle{(\frac{n}{2}-1)}d\varphi((\theta-\varphi)^{\sharp},Z)-d\varphi(J(\theta-\varphi)^{\sharp},JZ)-(\theta-\varphi)(JZ)d\varphi(J^{\wedge})\\[6pt]
\hfill-\rho_{D}((\theta-\varphi)^{\sharp},Z)+\rho^{\ast}_{D}(J(\theta-\varphi)^{\sharp},JZ)=0,

\end{array}
$$

where $\sharp:T^{\ast}M\to TM$ is the isomorphism defined by means
the metric $g$.

\end{theorem}

\begin{cor}\label{dim 4}
In the notation of Theorem~\ref{psd-harm-int}, if $dim\,M=4$, $J$
determines a $(D,D')$-pseudo-harmonic map from $M$ into the positive
twistor space ${\mathcal Z}$ of $(M,{\mathfrak c})$ if and only if
the $2$-form $d(\theta-\varphi)$ is of type $(1,1)$ with respect to
$J$ and
$$
(\theta-\varphi)(JZ)d\varphi(J^{\wedge})
+\rho_{D}((\theta-\varphi)^{\sharp},Z)-\rho^{\ast}_{D}(J(\theta-\varphi)^{\sharp},JZ)=0\quad
{\rm{for}}~~  Z\in TM.
$$
\end{cor}

\smallskip

\noindent {\bf Remark}.  If ${\mathfrak c}$ consists of a single
metric $g$ and we set $\varphi_{g}=0$, then $D$ is the Levi-Civita
connection of $g$ and $D'$ is the Levi-Civita connection of
$\widetilde g_t$. In this case, Corollary~\ref{dim 4} coincides with
\cite[Theorem 1]{DHM}.

\medskip

If $ D$ is the Weyl connection determined by $g$ and
$\varphi=\theta$, then conditions $(i)$ and $(ii)$ of
Theorem~\ref{psd-harm-int} are trivially satisfied, hence
$J:M\to{\mathcal Z}$ is a pseudo-harmonic map. In fact, this is a
consequence of Proposition~\ref{covder-dif} since ${ D}J=0$ by
(\ref{Weyl}) and (\ref{nablaJ}). The latter identities implies that
$D$ is the unique Weyl connection preserving the conformal class of
$g$ and such that $J$ is parallel (\cite{V76}).

\smallskip

 \noindent {\bf Example 1}.  It has been observed in
\cite{Tr} that every Inoue surface $M$ of type $S^0$ admits a
locally conformally K\"ahler metric $g$  for which the Lee form
$\theta$ is nowhere vanishing. %(see also \cite{DO}).
Define the metric $\widetilde g_t$ on the twistor space  by means of
the Levi-Civita connection of the metric $g$. It is shown in
\cite{DHM} that in this case the map $J: (M,g)\to ({\mathcal
Z},\widetilde g_t)$ is not harmonic. But, $J$ is pseudo-harmonic if
we use the Weyl connection $D$ defined by means of  $g$ and
$\varphi_{g}=\theta$.

\smallskip

Now,  recall  the construction of the Inoue surfaces of type $S^0$
(\cite{Inoue}). Let $A\in SL(3,\mathbb{Z})$ be a matrix with a real
eigenvalue $\alpha > 1$ and two complex eigenvalues $\beta$ and
$\overline{\beta}$, $\beta\neq\overline{\beta}$. Choose eigenvectors
$(a_1,a_2,a_3)\in{\mathbb R}^3$ and $(b_1,b_2,b_3)\in {\mathbb C}^3$
of $A$ corresponding to $\alpha$ and $\beta$, respectively. Then the
vectors $(a_1,a_2,a_3), (b_1,b_2,b_3),
(\overline{b_1},\overline{b_2},\overline{b_3})$  are
$\mathbb{C}$-linearly independent. Denote the upper-half plane in
$\mathbb{C}$ by ${\bf H}$ and let $\Gamma$ be the group of
holomorphic automorphisms of ${\bf H}\times {\mathbb C}$ generated
by
$$g_o:(w,z)\to (\alpha w,\beta z), \quad g_i:(w,z)\to (w+a_i,z+b_i), \>i=1,2,3 .$$
The group $\Gamma$ acts on ${\bf H}\times{\mathbb C}$ freely and
properly discontinuously.  Then $M=({\bf H}\times {\mathbb
C})/\Gamma$ is a compact complex surface known as Inoue surface of
type $S^0$.

Following \cite{Tr}, consider on ${\bf H}\times {\mathbb C}$ the
Hermitian metric
$$
g=\frac{1}{v^2}(du\otimes du+dv\otimes dv)+v(dx\otimes dx+dy\otimes
dy),\quad u+iv\in{\bf H}, \quad x+iy\in {\mathbb C}.
$$
This metric is invariant under the action of the group $\Gamma$, so
it descends to a Hermitian metric on $M$ which we denote again by
$g$. Instead on $M$, we  work with $\Gamma$-invariant objects on
${\bf H}\times{\mathbb C}$.  Let $\Omega$ be the fundamental
$2$-form of the Hermitian structure $(g,J)$ on ${\bf
H}\times{\mathbb C}$, $J$ being the standard complex structure. Then
$$ d\Omega=\frac{1}{v}dv\wedge\Omega. $$ Hence the Lee form is
$\theta=d\ln v$. In particular, $d\theta=0$, i.e. $(g,J)$ is a
locally conformally K\"ahler structure. Set
$$
E_1=v\frac{\partial}{\partial u},\quad E_2=v\frac{\partial}{\partial
v},\quad E_3=\frac{1}{\sqrt v}\frac{\partial}{\partial x},\quad
E_4=\frac{1}{\sqrt v}\frac{\partial}{\partial y}.
$$
These are $\Gamma$-invariant vector fields constituting an
orthonormal basis such that $JE_1=E_2$, $JE_3=E_4$. The non-zero Lie
brackets of $E_1,...,E_4$ are
$$
[E_1,E_2]=-E_1,\quad [E_2,E_3]=-\frac{1}{2}E_3,\quad
[E_2,E_4]=-\frac{1}{2}E_4.
$$

Denote  the dual basis of $E_1,...,E_4$ by $\eta_1,...,\eta_4$. The
$1$-forms $\eta_1,...,\eta_4$  are $\Gamma$-invariant. Take a
$\Gamma$-invariant $1$-form $\varphi=a_1\eta_1+...+a_4\eta_4$, where
$a_1,...,a_4$ are real constants.  The Levi-Civita connection
$\nabla$ of the metric $g$ is $\Gamma$-invariant, hence the
connection $D$ on ${\bf H}\times{\mathbb C}$ defined via
(\ref{Weyl})  by means the metric $g$ and the form $\varphi$ yields
a Weyl connection on $M$.

We have   $d\eta_1=\eta_1\wedge\eta_2$, $d\eta_2=0$,
$d\eta_3=\frac{1}{2}\eta_2\wedge \eta_3$,
$d\eta_4=\frac{1}{2}\eta_2\wedge \eta_4$ (the definition of the
exterior product being  so that $\eta_1\wedge \eta_2(E_1\wedge
E_2)=1$, etc.). Note also that $\theta=\eta_2$. Then the form
$$
d(\varphi-\theta)=a_1\eta_1\wedge\eta_2+\frac{1}{2}a_3\eta_2\wedge
\eta_3+\frac{1}{2}a_4\eta_2\wedge \eta_4
$$
is of type $(1,1)$ w.r.t. $J$ exactly when $a_3=a_4=0$. Thus $
\varphi=a_1\eta_1+a_2\eta_2$ and $d\varphi=a_1\eta_1\wedge\eta_2$.

The Levi-Civita connection $\nabla$ of $g$ is given by the following
table (\cite{DHM}):
\begin{equation}\label{LC-S0}
\begin{array}{c}
\nabla_{E_1}E_1=E_2,\quad \nabla_{E_1}E_2=-E_1,\\[6pt]
\displaystyle{\nabla_{E_3}E_2=\frac{1}{2}E_3,\quad
\nabla_{E_3}E_3=-\frac{1}{2}E_2, \quad
\nabla_{E_4}E_2=\frac{1}{2}E_4,\quad
\nabla_{E_4}E_4=-\frac{1}{2}E_2},\\[6pt]
{\rm{all~ other}}~ \nabla_{E_i}E_j=0.
\end{array}
\end{equation}
Now, since  $d\varphi=a_1\eta_1\wedge\eta_2$ is of type $(1,1)$ with
respect to $J$, condition $(ii)$ of Theorem~\ref{psd-harm-int} takes
the form
\begin{equation}\label{Ino}
(\theta-\varphi)(JZ)d\varphi(J^{\wedge})+\rho_{D}((\theta-\varphi)^{\sharp},Z)-\rho^{\ast}_{D}(J(\theta-\varphi)^{\sharp},JZ)=0.
\end{equation}
We have $d\varphi(J^{\wedge})=d\varphi(E_1\wedge
E_2)+d\varphi(E_3\wedge E_4)=a_1$.  Using (\ref{LC-S0}), we can
compute the curvature of the Weyl connection by means of identity
(\ref{RD-Rg}). As a result,  the non-zero components of
${\rho_D}_{ij}={\rho_D}(E_i,E_j)$ are
$$
\begin{array}{c}
\displaystyle{{\rho_D}_{11}=-\frac{a_2(a_2+2)}{2}, \quad {\rho_D}_{12}=\frac{a_1(a_2+3)}{2}, \quad {\rho_D}_{21}=\frac{a_1a_2}{2}}, \\[6pt]

\displaystyle{{\rho_D}_{22}=-\frac{a^2_1+3}{2}, \quad
{\rho_D}_{33}={\rho_D}_{44}=-\frac{a^2_1+a_2(a_2-1)}{2}}.
\end{array}
$$
Also, the non-zero components of the $\ast$-Ricci tensor of $(g,J)$
${\rho^{\ast}_D}_{ij}={\rho^{\ast}_D}(E_i,E_j)$ are the following:
$$
\begin{array}{c}
\displaystyle{{\rho^{\ast}_D}_{11}={\rho^{\ast}_D}_{22}=-\frac{a_2+2}{2}, \quad {\rho^{\ast}_D}_{12}=-{\rho^{\ast}_D}_{21}=-{\rho^{\ast}_D}_{34}={\rho^{\ast}_D}_{43}=\frac{a_1}{2}},\\[6pt]

\displaystyle{{\rho^{\ast}_D}_{33}={\rho^{\ast}_D}_{44}=-\frac{a^2_1+(a_2-1)^2}{4}}.
\end{array}
$$
Now, it is easy to see that identity (\ref{Ino}) holds for every $Z$
if and only if
$$
a_1(a_2-1)=0, \quad a_1^2+(a_2-1)^2=0.
$$
Obviously, the solution of the latter system is $a_1=0$, $a_2=1$.
Thus, $J$ is a pseudo-harmonic map from $M$ into its twistor space
only when $\varphi=\theta$, i.e. $D$ is the Weyl connection
determined by $g$ and the Lee form $\theta$ of $(g,J)$.

\smallskip

\noindent {\bf Example 2}. Recall that a primary Kodaira surface $M$
is the quotient of ${\mathbb C}^2$ by  a group of transformations
acting freely and properly discontinuously \cite[p. 787]{Kodaira}.
This group is generated by the affine transformations
$\varphi_k(z,w)=(z+a_k,w+\overline{a}_kz+b_k)$, where $a_k$, $b _k$,
$k=1,2,3,4$, are complex numbers such that $ a_1=a_2=0$, $b_2\neq
0$, $Im(a_3{\overline a}_4)=m b_1\neq 0$ for some integer $m>0$. The
quotient space is compact.

It is well-known that $M$ can also be described  as the quotient of
${\mathbb C}^2$ endowed with a group structure by a discrete
subgroup $\Gamma$. The multiplication on ${\Bbb C}^2$ is defined by
$$
(a,b).(z,w)=(z+a,w+\overline{a}z+b),\quad (a,b), (z,w)\in  {\Bbb
C}^2,
$$
and $\Gamma$ is the subgroup generated by $(a_k,b_k)$, $k=1,...,4$
(see, for example, \cite{Borc}). Considering $M$ as the quotient
${\Bbb C}^2/{\Gamma}$, every left-invariant object on ${\Bbb C}^2$
descends to a globally defined object on $M$.

As in \cite{D14, DHM, Dav}, take a frame of left-invariant  vector
fields $A_1,...,A_4$ such that
$$
[A_1,A_2]=-2A_4,\quad [A_i,A_j]=0 ~~\rm{otherwise.}
$$
These identities are satisfied, for example, by the following
left-invariant frame
\begin{equation}\label{A-frame}
A_1=-\frac{\partial}{\partial x}-x\frac{\partial}{\partial
u}+y\frac{\partial}{\partial v},\quad A_2=\frac{\partial}{\partial
y}+y\frac{\partial}{\partial u}+x\frac{\partial}{\partial v},\quad
A_3=\frac{\partial}{\partial u},\quad A_4=\frac{\partial}{\partial
v},
\end{equation}
where $x+iy=z$, $u+iv=w$.

Denote by $g$ the left-invariant Riemannian metric on $M$ for which
the frame $A_1,...,A_4$ is orthonormal.

By \cite{Has}, every complex structure on $M$ is induced by a
left-invariant complex structure on ${\Bbb C}^2$. It is easy to see
that every such a structure is given by (\cite{ M, D14})
$$
JA_1=\varepsilon_1 A_2,\quad JA_3=\varepsilon_2 A_4,\quad
\varepsilon_1,\varepsilon_2=\pm 1.
$$
Denote the complex structure defined by these identity by
$J_{\varepsilon_1,\varepsilon_2}$.

The non-zero covariant derivatives $\nabla_{A_i}A_j$ are
(\cite{DHM})
$$
\nabla_{A_1}A_2=-\nabla_{A_2}A_1=-A_4,\quad
\nabla_{A_1}A_4=\nabla_{A_4}A_1=A_2,\quad
\nabla_{A_2}A_4=\nabla_{A_4}A_2=-A_1.
$$
This implies that the Lee form is
$$
\theta(X)=-2\varepsilon_1g(X,A_3).
$$
Therefore $B=-2\varepsilon_1A_3,\quad \nabla\theta=0$.

Denote  the dual basis of $A_1,...,A_4$ by $\alpha_1,...,\alpha_4$.
Thus, $\theta=-2\varepsilon_1\alpha_3$.  Any left-invariant $1$-form
$\varphi$ is of the form $ \varphi=a_1\alpha_1+...+a_4\alpha_4$,
where $a_1,...,a_4$ are real constants. The connection $D$ on
${\mathbb C}^2$ defined via (\ref{Weyl})  by means the metric $g$
and the form $\varphi$ yields a Weyl connection on $M$.

We have $d\alpha_1=d\alpha_2=d\alpha_3=0$,
$d\alpha_4=2\alpha_1\wedge \alpha_2$. Hence the form $
d(\varphi-\theta)=d\varphi=2a_4\alpha_1\wedge \alpha_2$ is of type
$(1,1)$ w.r.t. $J=J_{\varepsilon_1,\varepsilon_2}$.

Let $R$ be the curvature tensor of the Levi-Civita connection
$\nabla$. Set for short $R_{ijk}=R(A_i,A_j)A_k$. Then the non-zero
$R_{ijk}$ are (\cite{DHM})
\begin{equation}\label{R-LC}
\begin{array}{c}
R_{121}=-3A_2,\quad R_{122}=3A_1,\quad R_{141}=A_4,\\[6pt]
 R_{144}=-A_1,\quad R_{242}=A_4,\quad R_{244}=-A_2.
\end{array}
\end{equation}
The non-zero covariant derivatives of the form  $
\varphi=a_1\alpha_1+...+a_4\alpha_4$ are
$$
\begin{array}{c}
(\nabla_{A_1}\varphi)(A_2)=-(\nabla_{A_2}\varphi)(A_1)=a_4, \quad
(\nabla_{A_1}\varphi)(A_4)=(\nabla_{A_4}\varphi)(A_1)=-a_2,\\[6pt]
(\nabla_{A_2}\varphi)(A_4)=\nabla_{A_4}\varphi)(A_2)=a_1.
\end{array}
$$
Now, using (\ref{R-LC}), we can compute the curvature of the Weyl
connection $D$ by means of identity (\ref{RD-Rg}), then the
components ${\rho_D}_{ij}={\rho_D}(A_i,A_j)$ and
${\rho^{\ast}_D}_{ij}={\rho^{\ast}_D}(A_i,A_j)$ of the Ricci and
$\ast$-Ricci tensors. The latter are given in the following tables.
$$
\begin{array}{c}

\displaystyle{{\rho_D}_{11}=-\frac{a_2^2+a_3^2+a_4^2+4}{2}, \quad {\rho_D}_{12}=2a_4+\frac{1}{2}a_1a_2, \quad {\rho_D}_{13}={\rho_D}_{31}=\frac{1}{2}a_1a_3},\\[6pt]

\displaystyle{{\rho_D}_{14}=-a_2+\frac{1}{2}a_1a_4, \quad {\rho_D}_{21}=-2a_4+\frac{1}{2}a_1a_2, \quad {\rho_D}_{22}=-\frac{a_1^2+a_3^2+a_4^2+4}{2}},\\[6pt]

\displaystyle{{\rho_D}_{23}={\rho_D}_{32}=\frac{1}{2}a_2a_3, \quad {\rho_D}_{24}={\rho_D}_{42}=a_1+\frac{1}{2}a_2a_4, \quad {\rho_D}_{33}=-\frac{a_1^2+a_2^2+a_4^2}{2}},\\[6pt]

\displaystyle{{\rho_D}_{34}={\rho_D}_{43}=\frac{1}{2}a_3a_4, \quad
{\rho_D}_{41}=-a_2+\frac{1}{2}a_1a_4, \quad
{\rho_D}_{44}=-\frac{a_1^2+a_2^2+a_3^2-4}{2}}.

\end{array}
$$
$$
\begin{array}{c}

\displaystyle{{\rho^{\ast}_D}_{11}={\rho^{\ast}_D}_{22}=-\frac{a_3^2+a_4^2+12}{4}, \quad {\rho^{\ast}_D}_{12}=-{\rho^{\ast}_D}_{21}=a_4},\\[6pt]

\displaystyle{{\rho^{\ast}_D}_{13}={\rho^{\ast}_D}_{31}=\varepsilon_1\varepsilon_2{\rho^{\ast}_D}_{24}
=\varepsilon_1\varepsilon_2{\rho^{\ast}_D}_{42}=\frac{a_1a_3+\varepsilon_1\varepsilon_2(2a_1+a_2a_4)}{4}},\\[6pt]

\displaystyle{{\rho^{\ast}_D}_{14}={\rho^{\ast}_D}_{41}=-\varepsilon_1\varepsilon_2{\rho^{\ast}_D}_{23}
=-\varepsilon_1\varepsilon_2{\rho^{\ast}_D}_{32}=\frac{-2a_2+a_1a_4-\varepsilon_1\varepsilon_2a_2a_3}{4}},\\[6pt]

\displaystyle{{\rho^{\ast}_D}_{34}=-{\rho^{\ast}_D}_{43}=-\varepsilon_1\varepsilon_2a_4,
\quad
{\rho^{\ast}_D}_{33}={\rho^{\ast}_D}_{44}=-\frac{a_1^2+a_2^2}{4}}.

\end{array}
$$

Condition $(ii)$ of Theorem~\ref{psd-harm-int} reduces to the
identity
$$
\begin{array}{c}
2a_4\varepsilon_1(\theta-\varphi)(JZ)
+\rho_{D}((\theta-\varphi)^{\sharp},Z)-\rho^{\ast}_{D}(J(\theta-\varphi)^{\sharp},JZ)=0.
\end{array}
$$
This identity is satisfied for every $Z\in TM$ if and only if

$$
\begin{array}{c}

(1-\varepsilon_2)a_2a_4+(1-\varepsilon_2)(2+\varepsilon_1 a_3)a_1=0,\\[6pt]

(1-\varepsilon_2)a_1a_4-(1-\varepsilon_2)(2+\varepsilon_1a_3)a_2=0,\\[6pt]

(1+\varepsilon_2)(a_1^2+a_2^2)+ 2(1-\varepsilon_2)a_4^2=0,\\[6pt]

(1-\varepsilon_2)(2+\varepsilon_1a_3)a_4=0.
\end{array}
$$
Clearly, if  $\varepsilon_2=1$, the solution of this system is
$a_1=a_2=0$, $a_3$ and $a_4$ arbitrary. For $\varepsilon_2=-1$, the
solutions are  $a_1$, $a_2$ arbitrary, $a_3=-2\varepsilon_1$,
$a_4=0$ and $a_1=a_2=a_4=0$, $a_3\neq -2\varepsilon_1$.  Thus, on a
Kodaira surface, there are many Weyl connections for which the
complex structures $J_{\varepsilon_1,\varepsilon_2}$ are
pseudo-harmonic maps from $M$ into its twistor space ${\mathcal Z}$.

\medskip

The next statement is a weaker version of \cite[Proposition
1.3]{Pon}.

\begin{cor}\label{unique}
Let $J_1$ and $J_2$ be two  Hermitian structures on a (connected) Riemannian manifold $(M,g)$.
Suppose that $J_1$ and $J_2$ have the same Lee form $\theta$ satisfying identity (\ref{dOm}). If $J_1$ and $J_2$ coincide on an open subset or, more-generally, if they coincide to infinite order at a point, they coincide on the whole manifold $M$.
\end{cor}

\begin{proof} The complex structures $J_1$ and $J_2$ determine the same orientation on $M$,
hence the same positive twistor space ${\mathcal Z}$. Let  $ D$ be
the Weyl connection determined $g$ and $\theta$. If $\widetilde g_t$
is the metric on ${\mathcal Z}$ defined by means of $ D$, and let
$D'$ be the torsion-free connection on ${\mathcal Z}$ defined by
(\ref{new conn}). Then the maps $J_1$ and $J_2$ of $M$ into
${\mathcal Z}$ are $(D,D')$-pseudo-harmonic by
Theorem~\ref{psd-harm-int}. Thus, the result follows from the
uniqueness theorem for pseudo-harmonic maps \cite{K09}.
\end{proof}

\smallskip

\noindent {\bf Remark}. In fact, M. Pontecorvo \cite{Pon} has proved
that the conclusion of the above corollary holds without the
assumption on the Lee forms. The idea of his proof is inspired by
the theory of pseudo-holomorphic curves. It makes use of the
observation in \cite{ES} that an almost-Hermitian structure $J$ on a
Riemannian manifold $M$ is integrable if and only if the
corresponding map $J: M\to {\mathcal Z}$ is pseudo-holomorphic with
respect to the almost-complex structure $J$ on $M$ and the
Atiyah-Hitchin-Singer \cite {AHS} almost-complex structure on
${\mathcal Z}$.

\bigskip

\noindent {\bf Acknowledgements}

\medskip

The authors would like to thank very much the referee, whose remarks
helped to correct and improve the final version of the paper.

\smallskip

The second named author is partially supported by the Bulgarian
National Science Fund, Ministry of Education and Science of Bulgaria
under contract DN 12/2.

\smallskip

This paper has been completed during the visit of the second name
author in the Abdus Salam School of Mathematical Sciences,
Government College University, Lahore, Pakistan. He would like to
thank  the whole staff of the School for their hospitality.

\end{document}